\documentclass[12pt,leqno]{article}

\usepackage{theorem,amssymb,amstex,amsopn}

\makeatletter
\ifx\@auxdefsloaded\relax\endinput
\else\let\@auxdefsloaded\relax\fi

\def\newenvironment{%
   \@ifnextchar *{\@@newenv{\global\@ignoretrue}}{\@@newenv{}*}}

\def\@@newenv#1*#2{%
   \@ifnextchar [{\@newenv{#1}{#2}}{\@newenv{#1}{#2}[0]}}

\long\def\@newenv#1#2[#3]#4#5{%
   \expandafter\newcommand\csname#2\endcsname[#3]{#4}%
   \expandafter\long\expandafter\def\csname end#2\endcsname{#5#1}}

\def\renewenvironment{%
   \@ifnextchar *{\@@renewenv{\global\@ignoretrue}}{\@@renewenv{}*}}

\def\@@renewenv#1*#2{%
   \@ifnextchar [{\@renewenv{#1}{#2}}{\@renewenv{#1}{#2}[0]}}

\long\def\@renewenv#1#2[#3]#4#5{%
   \expandafter\renewcommand\csname#2\endcsname[#3]{#4}%
   \expandafter\long\expandafter\def\csname end#2\endcsname{#5#1}}

\def\newoptcommand#1#2{%
   \@ifnextchar [{\@optargdef#1#2}{\@optargdef#1#2[1]}}

\def\renewoptcommand#1#2{%
   \edef\@tempa{\expandafter\@cdr\string#1\@nil}%
   \@ifundefined{\@tempa}{%
      \@latexerr{\string#1\space undefined}\@ehc}{}%
   \@ifnextchar [{\@reoptargdef#1#2}{\@reoptargdef#1#2[1]}}

\long\def\@optargdef#1#2[#3]#4{%
   \@ifdefinable #1{\@reoptargdef#1#2[#3]{#4}}}

\long\def\@reoptargdef#1#2[#3]#4{%
   \@tempcnta#3\relax \@tempcntb \@ne
   \let#1\relax \let\@tempa\relax
   \edef\@tempb{\long\def\csname\string#1\endcsname
      [\@tempa\the\@tempcntb]}%
   \advance\@tempcntb \@ne \advance\@tempcnta \m@ne
   \@whilenum\@tempcnta>0\do{%
      \edef\@tempb{\@tempb\@tempa\the\@tempcntb}%
      \advance\@tempcntb \@ne \advance\@tempcnta \m@ne}%
   \let\@tempa=##\@tempb{#4}%
   \def#1{\@ifnextchar [{\csname\string#1\endcsname}{%
      \csname\string#1\endcsname[#2]}}}

\def\newoptenvironment{%
   \@ifnextchar *{\@@newoptenv{\global\@ignoretrue}}{%
      \@@newoptenv{}*}}

\def\@@newoptenv#1*#2#3{%
   \@ifnextchar [{\@newoptenv{#1}{#2}{#3}}{%
      \@newoptenv{#1}{#2}{#3}[0]}}

\long\def\@newoptenv#1#2#3[#4]#5#6{%
   \expandafter\newoptcommand\csname#2\endcsname{#3}[#4]{#5}%
   \expandafter\long\expandafter\def\csname end#2\endcsname{#6#1}}

\def\renewoptenvironment{%
   \@ifnextchar *{\@@renewoptenv{\global\@ignoretrue}}{%
      \@@renewoptenv{}*}}

\def\@@renewoptenv#1*#2#3{%
   \@ifnextchar [{\@renewoptenv{#1}{#2}{#3}}{%
      \@renewoptenv{#1}{#2}{#3}[0]}}

\long\def\@renewoptenv#1#2#3[#4]#5#6{%
   \expandafter\renewoptcommand\csname#2\endcsname{#3}[#4]{#5}%
   \expandafter\long\expandafter\def\csname end#2\endcsname{#6#1}}

\newcounter{keepoptional}
\newcounter{optnestctr}

\newoptenvironment*{optional}{1}[1]{%
   \ifnum#1>\value{keepoptional}\relax
      \setcounter{optnestctr}{0}\@powerup\expandafter\@endopt\fi
   \ignorespaces}{}

\def\@powerup{\catcode`\{=12 \catcode`\}=12 \catcode`\\=12 \relax}
\def\@powerdown{\catcode`\{=1 \catcode`\}=2 \catcode`\\=0 \relax}
\begingroup \catcode`|=0 \catcode`[=1 \catcode`]=2 \@powerup
   |long|gdef|@endopt#1\end{optional}[%
      |@beginopt#1\begin{optional}*]%
   |long|gdef|@beginopt#1\begin{optional}[%
      |@ifstar[%
         |ifnum|value[optnestctr]>0|relax
            |addtocounter[optnestctr][-1]|let|@zzzap=|@endopt
         |else
            |@powerdown|end[optional]|let|@zzzap=|relax|fi
         |@zzzap]%
      [
         |addtocounter[optnestctr][1]|@beginopt]]%
|endgroup

\makeatother
\makeatletter
\ifx\@auxdefsloaded\relax
\else \input{auxdefs.sty}\fi
\newskip\dgARROWLENGTH  \dgARROWLENGTH=2.5em\relax
\newskip\dgHORIZPAD     \dgHORIZPAD=1em\relax
\newskip\dgVERTPAD      \dgVERTPAD=2ex\relax
\newskip\dgLABELOFFSET  \dgLABELOFFSET=.7ex\relax
\newcount\dgARROWPARTS  \dgARROWPARTS=4\relax
\newcommand{\dgeverynode}{\displaystyle}
\newcommand{\dgeverylabel}{\scriptstyle}
\newskip\dgDOTSPACING   \dgDOTSPACING=0.35em
\newskip\dgDOTSIZE      \dgDOTSIZE=1.5\fontdimen8\tenln
\newcount\dgMAXSQUARE   \dgMAXSQUARE=4\relax
\newcount\dgMINDBLSQ    \dgMINDBLSQ=10\relax
\newskip\dgCOLUMNWIDTH  \dgCOLUMNWIDTH=2em\relax
\chardef\f@ur=4
\@namedef{dgo@}{\let\dg@VECTOR=\vector
   \dg@LBLPOS=\dgARROWPARTS \divide\dg@LBLPOS\tw@}%
\@namedef{dgo@-}{\let\dg@VECTOR=\line}%
\@namedef{dgo@!}{\let\dg@VECTOR=\dg@novector}%
\@namedef{dgo@1}{\dg@LBLPOS=\@ne\relax}%
\@namedef{dgo@2}{\dg@LBLPOS=\tw@\relax}%
\@namedef{dgo@3}{\dg@LBLPOS=\thr@@\relax}%
\@namedef{dgo@4}{\dg@LBLPOS=\f@ur\relax}%
\@namedef{dgo@5}{\dg@LBLPOS=5\relax}%
\@namedef{dgo@6}{\dg@LBLPOS=6\relax}%
\@namedef{dgo@7}{\dg@LBLPOS=7\relax}%
\@namedef{dgo@8}{\dg@LBLPOS=8\relax}%
\@namedef{dgo@9}{\dg@LBLPOS=9\relax}%
\def\dgt@e{\dg@DX=\@ne \dg@DY=\z@ \dg@SIZE=\@ne}%
\def\dgt@w{\dg@DX=\m@ne \dg@DY=\z@ \dg@SIZE=\@ne}%
\def\dgt@n{\dg@DX=\z@ \dg@DY=\@ne \dg@SIZE=\@ne}%
\def\dgt@s{\dg@DX=\z@ \dg@DY=\m@ne \dg@SIZE=\@ne}%
\def\dgt@ne{\dg@DX=\@ne \dg@DY=\@ne \dg@SIZE=\@ne}%
\def\dgt@se{\dg@DX=\@ne \dg@DY=\m@ne \dg@SIZE=\@ne}%
\def\dgt@nw{\dg@DX=\m@ne \dg@DY=\@ne \dg@SIZE=\@ne}%
\def\dgt@sw{\dg@DX=\m@ne \dg@DY=\m@ne \dg@SIZE=\@ne}%
\def\dgt@nne{\dg@DX=\@ne \dg@DY=\tw@ \dg@SIZE=\@ne}%
\def\dgt@nnw{\dg@DX=\m@ne \dg@DY=\tw@ \dg@SIZE=\@ne}%
\def\dgt@sse{\dg@DX=\@ne \dg@DY=-\tw@ \dg@SIZE=\@ne}%
\def\dgt@ssw{\dg@DX=\m@ne \dg@DY=-\tw@ \dg@SIZE=\@ne}%
\def\dgt@ene{\dg@DX=\tw@ \dg@DY=\@ne \dg@SIZE=\tw@}%
\def\dgt@ese{\dg@DX=\tw@ \dg@DY=\m@ne \dg@SIZE=\tw@}%
\def\dgt@wnw{\dg@DX=-\tw@ \dg@DY=\@ne \dg@SIZE=\tw@}%
\def\dgt@wsw{\dg@DX=-\tw@ \dg@DY=\m@ne \dg@SIZE=\tw@}%
\def\dggeometry{
   \dg@ZTEMP=\dg@XGRID \multiply\dg@ZTEMP\tw@
   \ifnum\dg@YGRID=\z@ \dg@ZTEMP=\tw@
   \else \divide\dg@ZTEMP\dg@YGRID \fi
   \ifnum\dg@ZTEMP>\f@ur \dg@ZTEMP=\f@ur \fi
   \ifnum\dg@ZTEMP<\@ne \dg@ZTEMP=\@ne \fi
   \unitlength=2sp\relax
   \ifnum\dg@ZTEMP<\tw@
      \advance\dg@ZTEMP\@ne
      \multiply\unitlength\dg@YGRID
   \else
      \multiply\unitlength\dg@XGRID \divide\unitlength\dg@ZTEMP
   \fi
   \dg@XGRID=\dg@ZTEMP \dg@YGRID=\tw@
   \dg@rmcommondiv\tw@\dg@XGRID\dg@YGRID
   \divide\unitlength\dg@YGRID \divide\unitlength\@m\relax}
\@namedef{dgo@..}{\let\dg@VECTOR=\dg@dotvector}%

\def\dg@dotvector(#1,#2)#3{%
   \begingroup
   \dg@XTEMP=#1\relax \dg@YTEMP=#2\relax
   \let\dg@NDOTS=\dg@XEND \let\dg@DOTDIAM=\dg@WEND
   \dg@NDOTS=\unitlength \multiply\dg@NDOTS #3\relax
   \dg@ZTEMP=\dg@YTEMP \dg@changesign\dg@YTEMP\dg@ZTEMP
   \ifnum\dg@XTEMP>\z@
      \ifnum\dg@YTEMP>\dg@XTEMP
         \multiply\dg@NDOTS\dg@YTEMP \divide\dg@NDOTS\dg@XTEMP \fi
   \else\ifnum\dg@XTEMP<\z@
      \ifnum\dg@YTEMP>-\dg@XTEMP
         \multiply\dg@NDOTS\dg@YTEMP \divide\dg@NDOTS-\dg@XTEMP \fi
   \fi\fi
   \dg@YTEMP=\dg@ZTEMP
   \divide\dg@NDOTS\dgDOTSPACING
   \ifnum\dg@NDOTS>\z@\else \dg@NDOTS=\@ne \fi
   \dg@ZTEMP=\unitlength \multiply\dg@ZTEMP #3\relax
   \divide\dg@ZTEMP\dg@NDOTS
   \ifnum\dg@XTEMP=\z@
      \dg@changesign\dg@ZTEMP\dg@YTEMP \dg@YTEMP=\dg@ZTEMP
   \else
      \dg@changesign\dg@ZTEMP\dg@XTEMP
      \multiply\dg@YTEMP\dg@ZTEMP \divide\dg@YTEMP\dg@XTEMP
      \dg@XTEMP=\dg@ZTEMP
   \fi
   \divide\dg@XTEMP\unitlength \divide\dg@YTEMP\unitlength
   \begin{picture}(0,0)
      \dg@DOTDIAM=\dgDOTSIZE \divide\dg@DOTDIAM\unitlength
      \multiput(0,0)(\dg@XTEMP,\dg@YTEMP){\dg@NDOTS}{%
         \circle*{\dg@DOTDIAM}}%
      \multiply\dg@XTEMP\dg@NDOTS \multiply\dg@YTEMP\dg@NDOTS
      \put(\dg@XTEMP,\dg@YTEMP){\vector(#1,#2){0}}%
   \end{picture}%
   \endgroup}%
\newif\ifdg@LATEXGEOM \dg@LATEXGEOMfalse

\@ifundefined{lamsvector}{}{%
   \@namedef{dgo@}{%
      \let\dg@VECTOR=\lamsvector
      \lamsshaft{?}\lamssource{?}\lamstarget{?}%
      \dg@LBLPOS=\dgARROWPARTS \divide\dg@LBLPOS\tw@\relax}%
   \@namedef{dgo@..}{\lamsshaft{.}}%
   \@namedef{dgo@=}{\lamsshaft{=}}%
   \@namedef{dgo@-}{\lamstarget{0}}%
   \@namedef{dgo@'}{\lamstarget{'}}%
   \@namedef{dgo@`}{\lamstarget{`}}%
   \@namedef{dgo@A}{\lamstarget{A}}%
   \@namedef{dgo@V}{\lamssource{V}}%
   \@namedef{dgo@J}{\lamssource{J}}%
   \@namedef{dgo@L}{\lamssource{L}}%
   \@namedef{dgo@S}{\lamssource{S}}%
   \def\dggeometry{
      \dg@ZTEMP=\dg@XGRID \multiply\dg@ZTEMP\tw@
      \ifnum\dg@YGRID=\z@ \dg@ZTEMP=\tw@
      \else \divide\dg@ZTEMP\dg@YGRID \fi
      \ifnum\dg@ZTEMP>6\relax \dg@ZTEMP=6\relax \fi
      \ifdg@LATEXGEOM\ifnum\dg@ZTEMP>\f@ur \dg@ZTEMP=\f@ur \fi\fi
      \ifnum\dg@ZTEMP<\@ne \dg@ZTEMP=\@ne \fi
      \unitlength=2sp\relax
      \ifnum\dg@ZTEMP<\tw@
         \advance\dg@ZTEMP\@ne
         \multiply\unitlength\dg@YGRID
      \else
         \multiply\unitlength\dg@XGRID \divide\unitlength\dg@ZTEMP
      \fi
      \dg@XGRID=\dg@ZTEMP \dg@YGRID=\tw@
      \dg@rmcommondiv\tw@\dg@XGRID\dg@YGRID
      \divide\unitlength\dg@YGRID \divide\unitlength\@m
      \dg@LATEXGEOMfalse}
   \def\dgt@nee{\dg@DX=\tw@ \dg@DY=\@ne \dg@SIZE=\tw@}%
   \def\dgt@see{\dg@DX=\tw@ \dg@DY=\m@ne \dg@SIZE=\tw@}%
   \def\dgt@nww{\dg@DX=-\tw@ \dg@DY=\@ne \dg@SIZE=\tw@}%
   \def\dgt@sww{\dg@DX=-\tw@ \dg@DY=\m@ne \dg@SIZE=\tw@}%
   \def\dgt@nnne{\dg@DX=\@ne \dg@DY=\thr@@ \dg@SIZE=\@ne}%
   \def\dgt@nnnw{\dg@DX=\m@ne \dg@DY=\thr@@ \dg@SIZE=\@ne}%
   \def\dgt@sssw{\dg@DX=\m@ne \dg@DY=-\thr@@ \dg@SIZE=\@ne}%
   \def\dgt@ssse{\dg@DX=\@ne \dg@DY=-\thr@@ \dg@SIZE=\@ne}%
   \def\dgt@nnnee{\dg@DX=\tw@ \dg@DY=\thr@@ \dg@SIZE=\tw@}%
   \def\dgt@nnnww{\dg@DX=-\tw@ \dg@DY=\thr@@ \dg@SIZE=\tw@}%
   \def\dgt@sssww{\dg@DX=-\tw@ \dg@DY=-\thr@@ \dg@SIZE=\tw@}%
   \def\dgt@sssee{\dg@DX=\tw@ \dg@DY=-\thr@@ \dg@SIZE=\tw@}%
   \def\dgt@nneee{\dg@DX=\thr@@ \dg@DY=\tw@ \dg@SIZE=\thr@@}%
   \def\dgt@nnwww{\dg@DX=-\thr@@ \dg@DY=\tw@ \dg@SIZE=\thr@@}%
   \def\dgt@sswww{\dg@DX=-\thr@@ \dg@DY=-\tw@ \dg@SIZE=\thr@@}%
   \def\dgt@sseee{\dg@DX=\thr@@ \dg@DY=-\tw@ \dg@SIZE=\thr@@}%
   \def\dgt@neee{\dg@DX=\thr@@ \dg@DY=\@ne \dg@SIZE=\thr@@
      \global\dg@LATEXGEOMtrue}%
   \def\dgt@nwww{\dg@DX=-\thr@@ \dg@DY=\@ne \dg@SIZE=\thr@@
      \global\dg@LATEXGEOMtrue}%
   \def\dgt@swww{\dg@DX=-\thr@@ \dg@DY=\m@ne \dg@SIZE=\thr@@
      \global\dg@LATEXGEOMtrue}%
   \def\dgt@seee{\dg@DX=\thr@@ \dg@DY=\m@ne \dg@SIZE=\thr@@
      \global\dg@LATEXGEOMtrue}%
}
\newcount\dg@HORIZ      \newcount\dg@VERT
\newcount\dg@XLPAD      \newcount\dg@YBPAD
\newcount\dg@XRPAD      \newcount\dg@YTPAD
\newcount\dg@X          \newcount\dg@Y
\newcount\dg@XGRID      \newcount\dg@YGRID
\newcount\dg@SIZE
\newcount\dg@USERSIZE
\newcount\dg@DX         \newcount\dg@DY
\newcount\dg@XLBL       \newcount\dg@YLBL
\newcount\dg@XOFFSET    \newcount\dg@YOFFSET
\newcount\dg@LBLOFF
\newcount\dg@XLBLOFF    \newcount\dg@YLBLOFF
\newcount\dg@LBLPOS
\newcount\dg@XTEMP      \newcount\dg@YTEMP
\newcount\dg@ZTEMP
\newcount\dg@XNODE      \newcount\dg@YNODE
\newcount\dg@XEND       \newcount\dg@YEND
\newcount\dg@WEND       \newcount\dg@HEND
\newcount\dg@COUNT
\newbox\dg@NODEBOX
\newoptenvironment*{diagram}{}[1]{%
   \global\advance\dg@COUNT\@ne \typeout{[diagram \the\dg@COUNT:}%
   \let\node=\dg@node \let\\=\dg@cr \let\arrow=\dg@arrow
   \def\dg@BIGNODE{#1}%
   \ifx\dg@BIGNODE\@empty\else
      \setbox\dg@NODEBOX=\hbox{$\dgeverynode{#1}$}%
      \dg@XTEMP=\wd\dg@NODEBOX
      \dg@YTEMP=\ht\dg@NODEBOX \advance\dg@YTEMP\dp\dg@NODEBOX
      \ifnum\dg@YTEMP=\z@ \dg@ZTEMP=\z@
      \else \dg@ZTEMP=\dg@XTEMP \divide\dg@ZTEMP\dg@YTEMP\fi
      \ifnum\dg@ZTEMP>\dgMAXSQUARE
         \ifnum\dg@ZTEMP<\dgMINDBLSQ \dg@XGRID=\thr@@ \dg@YGRID=\tw@
         \else \dg@XGRID=\tw@ \dg@YGRID=\@ne \fi
      \else \dg@XGRID=\@ne \dg@YGRID=\@ne \fi
      \advance\dg@XTEMP\dgHORIZPAD \advance\dg@YTEMP\dgVERTPAD
      \dg@XLPAD=\dg@XTEMP \divide\dg@XLPAD\tw@
      \dg@XRPAD=\dg@XTEMP \divide\dg@XRPAD\tw@
      \dg@YBPAD=\dg@YTEMP \divide\dg@YBPAD\tw@
      \dg@YTPAD=\dg@YTEMP \divide\dg@YTPAD\tw@
      \advance\dg@XTEMP\dgARROWLENGTH
      \divide\dg@XTEMP\dg@XGRID \divide\dg@XTEMP\@m
      \unitlength=1sp\relax \multiply\unitlength\dg@XTEMP
   \fi
   \typeout{saving...}%
   \dg@X=-\@ne \dg@Y=\z@ \dg@HORIZ=\z@\relax
   \let\dg@SLIST=\@empty
   \let\dg@NLIST=\@empty \let\dg@ALIST=\@empty
   \let\dg@PASS=\dg@savepass
}{
   \dg@VERT=-\dg@Y\relax
   \typeout{calculating...}%
   \ifx\dg@BIGNODE\@empty
      \dg@XGRID=\z@ \dg@YGRID=\z@
      \dg@XLPAD=\z@ \dg@XRPAD=\z@ \dg@YBPAD=\z@ \dg@YTPAD=\z@
      \let\dg@PASS=\dg@geompass
      \dg@NLIST \dg@ALIST
      \ifnum\dg@XGRID>\z@\else \dg@XGRID=\quad\relax \fi
      \ifnum\dg@YGRID>\z@\else \dg@YGRID=\quad\relax \fi
      \dggeometry
      \ifdim\unitlength>\z@\else \unitlength=\quad \fi
      \ifnum\dg@XGRID>\z@\else \dg@XGRID=\@ne \fi
      \ifnum\dg@YGRID>\z@\else \dg@YGRID=\@ne \fi
   \fi
   \dg@LBLOFF=\dgLABELOFFSET \divide\dg@LBLOFF\unitlength
   \multiply\dg@HORIZ\@m \multiply\dg@HORIZ\dg@XGRID
   \multiply\dg@VERT\@m \multiply\dg@VERT\dg@YGRID
   \divide\dg@XLPAD\unitlength \divide\dg@XRPAD\unitlength
   \divide\dg@YBPAD\unitlength \divide\dg@YTPAD\unitlength
   \typeout{drawing...}%
   \let\dg@PASS=\dg@drawpass
   \hspace*{\dg@XLPAD\unitlength}%
   \vcenter{%
      \hsize=0pt\relax\parindent=0pt\relax
      \parskip=0pt\relax\baselineskip=0pt\relax
      \vspace*{\dg@YTPAD\unitlength}%
      \begin{picture}(0,0)(0,0)\dg@NLIST\dg@ALIST\end{picture}%
      \raisebox{\z@}[\z@][\dg@VERT\unitlength]{}%
      \vspace*{\dg@YBPAD\unitlength}%
   }
   \hspace*{\dg@HORIZ\unitlength}%
   \hspace*{\dg@XRPAD\unitlength}%
   \typeout{done]}}%

\def\dg@savepass{s}
\def\dg@geompass{g}
\def\dg@drawpass{d}

\newoptcommand{\dg@node}{\@ne}[2]{%
   \ifx\dg@PASS\dg@savepass
      \dg@XTEMP=#1\relax
      \ifnum\dg@XTEMP<\@ne \dg@XTEMP=\@ne\fi
      \advance\dg@X\dg@XTEMP
      \ifnum\dg@HORIZ<\dg@X \dg@HORIZ=\dg@X \fi
      \setbox\dg@NODEBOX=\hbox{$\dgeverynode{#2}$}%
      \dg@XTEMP=\wd\dg@NODEBOX \advance\dg@XTEMP\dgHORIZPAD
      \dg@YTEMP=\ht\dg@NODEBOX \advance\dg@YTEMP\dp\dg@NODEBOX
      \advance\dg@YTEMP\dgVERTPAD
      \toks\z@=\expandafter{\dg@SLIST}%
      \edef\dg@SLIST{\the\toks\z@
         ,\noexpand\dg@XNODE=\number\dg@X\noexpand\relax
         \noexpand\dg@YNODE=\number\dg@Y\noexpand\relax
         \noexpand\dg@XTEMP=\number\dg@XTEMP\noexpand\relax
         \noexpand\dg@YTEMP=\number\dg@YTEMP\noexpand\relax}%
      \toks\z@=\expandafter{\dg@NLIST}%
      \toks\tw@={\dg@node{#2}}%
      \edef\dg@NLIST{\the\toks\z@
         \noexpand\dg@X=\number\dg@X\noexpand\relax
         \noexpand\dg@Y=\number\dg@Y\noexpand\relax
         \the\toks\tw@}%
   \else\ifx\dg@PASS\dg@geompass
      \ifnum\dg@X=\z@
         \dg@getnodesize
            {\dg@SLIST}{\dg@X}{\dg@Y}{\dg@WEND}{\dg@HEND}%
         \divide\dg@WEND\tw@
         \ifnum\dg@XLPAD<\dg@WEND \dg@XLPAD=\dg@WEND \fi\fi
      \ifnum\dg@X=\dg@HORIZ
         \dg@getnodesize
            {\dg@SLIST}{\dg@X}{\dg@Y}{\dg@WEND}{\dg@HEND}%
         \divide\dg@WEND\tw@
         \ifnum\dg@XRPAD<\dg@WEND \dg@XRPAD=\dg@WEND \fi\fi
      \ifnum\dg@Y=\z@
         \dg@getnodesize
            {\dg@SLIST}{\dg@X}{\dg@Y}{\dg@WEND}{\dg@HEND}%
         \divide\dg@HEND\tw@
         \ifnum\dg@YTPAD<\dg@HEND \dg@YTPAD=\dg@HEND \fi\fi
      \ifnum\dg@Y=-\dg@VERT\relax
         \dg@getnodesize
            {\dg@SLIST}{\dg@X}{\dg@Y}{\dg@WEND}{\dg@HEND}%
         \divide\dg@HEND\tw@
         \ifnum\dg@YBPAD<\dg@HEND \dg@YBPAD=\dg@HEND \fi\fi
   \else\ifx\dg@PASS\dg@drawpass
      \dg@XNODE=\dg@X \multiply\dg@XNODE\@m
      \multiply\dg@XNODE\dg@XGRID
      \dg@YNODE=\dg@Y \multiply\dg@YNODE\@m
      \multiply\dg@YNODE\dg@YGRID
      \setbox\dg@NODEBOX=\hbox{$\dgeverynode{#2}$}%
      \put(\dg@XNODE,\dg@YNODE){\dg@makebox{\box\dg@NODEBOX}}%
   \fi\fi\fi}%

\newoptcommand{\dg@cr}{\@ne}[1]{%
   \ifx\dg@PASS\dg@savepass
      \dg@YTEMP=#1\relax
      \ifnum\dg@YTEMP<\@ne \dg@YTEMP=\@ne \fi
      \advance\dg@Y -\dg@YTEMP\relax
      \dg@X=-\@ne\relax\fi}%
\newoptcommand{\dg@arrow}{\@ne}[2]{%
   \begingroup
   \dg@USERSIZE=#1\relax
   \ifnum\dg@USERSIZE<\@ne \dg@USERSIZE=\@ne \fi
   \dg@parse{#2}%
   \ifx\dg@PASS\dg@savepass
      \ifx\dg@label\dgl@b \let\dg@label=\dgl@t \fi
      \ifx\dg@label\dgl@r \let\dg@label=\dgl@l \fi
      \let\dg@process=\dg@save
   \else\ifx\dg@PASS\dg@geompass
      \let\dg@process=\dg@ignore
      \dg@geomcalc
   \else\ifx\dg@PASS\dg@drawpass
      \let\dg@process=\dg@draw
      \dg@drawcalc
   \fi\fi\fi
   \dg@label{\dg@process{#1}{#2}}}%

\newoptcommand{\arrow}{\@ne}[2]{%
   \dg@parse{#2}%
   \ifx\dg@label\dgl@b \let\dg@label=\dgl@t \fi
   \ifx\dg@label\dgl@r \let\dg@label=\dgl@l \fi
   \dg@label{\dg@textarrow{#1}{#2}}}%

\def\dg@textarrow#1#2#3#4{%
   \mathop{{\dgHORIZPAD=0pt\relax\dgVERTPAD=0pt\relax
      \begin{diagram}
         \node{}\arrow[#1]{#2}{#3}{#4}\node{}
      \end{diagram}}}}

\def\dg@parse#1{%
   \let\dg@label=\dgl@ \dgo@
   \let\dg@type=\@empty \let\dg@lbltype=\@empty
   \@for\dg@list:=#1\do{%
      \ifx\dg@type\@empty \let\dg@type=\dg@list
      \else\ifx\dg@lbltype\@empty \let\dg@lbltype=\dg@list
         \@ifundefined{dgo@\dg@list}{}{\@nameuse{dgo@\dg@list}}%
      \else
         \@ifundefined{dgo@\dg@list}{}{\@nameuse{dgo@\dg@list}}%
      \fi\fi}%
   \@ifundefined{dgt@\dg@type}{\dgt@e}{\@nameuse{dgt@\dg@type}}%
   \@ifundefined{dgl@\dg@lbltype}{}{%
      \dg@letname\dg@label{dgl@\dg@lbltype}}}

\def\dg@draw#1#2#3#4{%
   \put(\dg@X,\dg@Y){\dg@makebox{%
      \begin{picture}(0,0)%
         \thinlines
         \put(\dg@XOFFSET,\dg@YOFFSET){%
            \dg@VECTOR(\dg@DX,\dg@DY){\dg@SIZE}}%
         \put(\dg@XLBL,\dg@YLBL){\dg@makebox{%
            \begin{picture}(0,0)%
               \put(\dg@XLBLOFF,\dg@YLBLOFF){%
                  \dg@makebox[\dg@LBLONE]{$\dgeverylabel{#3}$}}%
               \put(-\dg@XLBLOFF,-\dg@YLBLOFF){%
                  \dg@makebox[\dg@LBLTWO]{$\dgeverylabel{#4}$}}%
            \end{picture}}}%
      \end{picture}}}%
   \endgroup}%

\def\dg@save#1#2#3#4{%
   \endgroup 
   \toks\z@=\expandafter{\dg@ALIST}%
   \toks\tw@={\dg@arrow[#1]{#2}{#3}{#4}}%
   \edef\dg@ALIST{\the\toks\z@%
      \noexpand\dg@X=\number\dg@X\noexpand\relax
      \noexpand\dg@Y=\number\dg@Y\noexpand\relax
      \the\toks\tw@}}%

\def\dg@ignore#1#2#3#4{\endgroup}

\def\dg@geomcalc{%
   \dg@XEND=\dg@SIZE \multiply\dg@XEND\dg@USERSIZE
   \ifnum\dg@DX=\z@
      \dg@YEND=\dg@XEND \dg@XEND=\z@
      \dg@changesign\dg@YEND\dg@DY
   \else
      \dg@changesign\dg@XEND\dg@DX \dg@YEND=\dg@XEND
      \multiply\dg@YEND\dg@DY \divide\dg@YEND\dg@DX
   \fi
   \advance\dg@XEND\dg@X \advance\dg@YEND\dg@Y
   \dg@getnodesize
      {\dg@SLIST}{\dg@XEND}{\dg@YEND}{\dg@WEND}{\dg@HEND}%
   \dg@XOFFSET=\dg@WEND \dg@YOFFSET=\dg@HEND
   \dg@getnodesize
      {\dg@SLIST}{\dg@X}{\dg@Y}{\dg@WEND}{\dg@HEND}%
   \advance\dg@XOFFSET\dg@WEND \divide\dg@XOFFSET\tw@
   \advance\dg@YOFFSET\dg@HEND \divide\dg@YOFFSET\tw@
   \dg@XTEMP=\dgARROWLENGTH \dg@YTEMP=\dgARROWLENGTH
   \ifnum\dg@DX>\z@
      \dg@ZTEMP=\dg@DX \multiply\dg@XTEMP\dg@DX
   \else \dg@ZTEMP=-\dg@DX \multiply\dg@XTEMP -\dg@DX \fi
   \ifnum\dg@DY>\z@
      \advance\dg@ZTEMP\dg@DY \multiply\dg@YTEMP\dg@DY
   \else \advance\dg@ZTEMP -\dg@DY \multiply\dg@YTEMP -\dg@DY\fi
   \ifnum\dg@ZTEMP=\z@\else
      \divide\dg@XTEMP\dg@ZTEMP \divide\dg@YTEMP\dg@ZTEMP
      \advance\dg@XOFFSET\dg@XTEMP \advance\dg@YOFFSET\dg@YTEMP
   \fi
   \divide\dg@XOFFSET\dg@SIZE \divide\dg@XOFFSET\dg@USERSIZE
   \divide\dg@YOFFSET\dg@SIZE \divide\dg@YOFFSET\dg@USERSIZE
   \ifnum\dg@DX=\z@ \dg@XOFFSET=\z@ \fi
   \ifnum\dg@DY=\z@ \dg@YOFFSET=\z@ \fi
   \ifnum\dg@XGRID<\dg@XOFFSET \global\dg@XGRID=\dg@XOFFSET\fi
   \ifnum\dg@YGRID<\dg@YOFFSET \global\dg@YGRID=\dg@YOFFSET\fi
   \relax}

\def\dg@drawcalc{%
   \dg@XEND=\dg@SIZE \multiply\dg@XEND\dg@USERSIZE
   \ifnum\dg@DX=\z@
      \dg@YEND=\dg@XEND \dg@XEND=\z@
      \dg@changesign\dg@YEND\dg@DY
   \else
      \dg@changesign\dg@XEND\dg@DX \dg@YEND=\dg@XEND
      \multiply\dg@YEND\dg@DY \divide\dg@YEND\dg@DX
   \fi
   \advance\dg@XEND\dg@X \advance\dg@YEND\dg@Y
   \dg@getnodesize
      {\dg@SLIST}{\dg@XEND}{\dg@YEND}{\dg@WEND}{\dg@HEND}%
   \divide\dg@WEND\unitlength \divide\dg@HEND\unitlength
   \multiply\dg@DX\dg@XGRID \multiply\dg@DY\dg@YGRID
   \dg@rmcommondiv\tw@\dg@DX\dg@DY
   \dg@rmcommondiv\tw@\dg@DX\dg@DY 
   \dg@rmcommondiv\thr@@\dg@DX\dg@DY
   \multiply\dg@SIZE\dg@USERSIZE \multiply\dg@SIZE\@m
   \ifnum\dg@DX=\z@
      \multiply\dg@SIZE\dg@YGRID
      \divide\dg@HEND\tw@ \advance\dg@SIZE -\dg@HEND
      \dg@getnodesize
         {\dg@SLIST}{\dg@X}{\dg@Y}{\dg@WEND}{\dg@YOFFSET}%
      \divide\dg@YOFFSET\unitlength \divide\dg@YOFFSET\tw@
      \advance\dg@SIZE -\dg@YOFFSET
      \dg@XOFFSET=\z@
      \def\dg@LBLONE{r}\def\dg@LBLTWO{l}%
      \dg@XLBL=\z@ \dg@YLBL=\dg@SIZE
      \multiply\dg@YLBL\dg@LBLPOS
      \divide\dg@YLBL\dgARROWPARTS\relax
      \advance\dg@YLBL\dg@YOFFSET
      \dg@changesign\dg@YLBL\dg@DY
      \dg@changesign\dg@YOFFSET\dg@DY
   \else
      \multiply\dg@SIZE\dg@XGRID
      \ifnum\dg@DY=\z@
         \divide\dg@WEND\tw@ \advance\dg@SIZE -\dg@WEND
         \dg@getnodesize
            {\dg@SLIST}{\dg@X}{\dg@Y}{\dg@XOFFSET}{\dg@HEND}%
         \divide\dg@XOFFSET\unitlength \divide\dg@XOFFSET\tw@
         \advance\dg@SIZE -\dg@XOFFSET
         \dg@YOFFSET=\z@
         \def\dg@LBLONE{b}\def\dg@LBLTWO{t}%
         \dg@YLBL=\z@ \dg@XLBL=\dg@SIZE
         \multiply\dg@XLBL\dg@LBLPOS
         \divide\dg@XLBL\dgARROWPARTS\relax
         \advance\dg@XLBL\dg@XOFFSET
         \dg@changesign\dg@XLBL\dg@DX
         \dg@changesign\dg@XOFFSET\dg@DX
      \else
         \divide\dg@WEND\tw@ \divide\dg@HEND\tw@
         \multiply\dg@HEND\dg@DX \divide\dg@HEND\dg@DY
         \ifnum\dg@HEND<\z@ \multiply\dg@HEND\m@ne \fi
         \ifnum\dg@WEND<\dg@HEND \advance\dg@SIZE -\dg@WEND
         \else \advance\dg@SIZE -\dg@HEND \fi
         \dg@getnodesize
            {\dg@SLIST}{\dg@X}{\dg@Y}{\dg@WEND}{\dg@HEND}%
         \divide\dg@WEND\unitlength \divide\dg@WEND\tw@
         \divide\dg@HEND\unitlength \divide\dg@HEND\tw@
         \multiply\dg@HEND\dg@DX \divide\dg@HEND\dg@DY
         \ifnum\dg@HEND<\z@ \multiply\dg@HEND\m@ne \fi
         \ifnum\dg@WEND<\dg@HEND \dg@XOFFSET=\dg@WEND
         \else \dg@XOFFSET=\dg@HEND \fi
         \advance\dg@SIZE -\dg@XOFFSET
         \dg@changesign\dg@XOFFSET\dg@DX
         \dg@YOFFSET=\dg@XOFFSET
         \multiply\dg@YOFFSET\dg@DY \divide\dg@YOFFSET\dg@DX
         \def\dg@LBLONE{br}\def\dg@LBLTWO{tl}%
         \ifnum\dg@DX<\z@ \ifnum\dg@DY>\z@
            \def\dg@LBLONE{bl}\def\dg@LBLTWO{tr}\fi\fi
         \ifnum\dg@DX>\z@ \ifnum\dg@DY<\z@
            \def\dg@LBLONE{bl}\def\dg@LBLTWO{tr}\fi\fi
         \dg@XLBL=\dg@SIZE
         \multiply\dg@XLBL\dg@LBLPOS
         \divide\dg@XLBL\dgARROWPARTS\relax
         \dg@changesign\dg@XLBL\dg@DX
         \dg@YLBL=\dg@XLBL
         \multiply\dg@YLBL\dg@DY \divide\dg@YLBL\dg@DX
         \advance\dg@XLBL\dg@XOFFSET
         \advance\dg@YLBL\dg@YOFFSET
      \fi
   \fi
   \dg@XLBLOFF=-\dg@DY \dg@changesign\dg@XLBLOFF\dg@DX
   \dg@YLBLOFF=\dg@DX \dg@changesign\dg@YLBLOFF\dg@DX
   \ifnum\dg@DX=\z@ \dg@XLBLOFF=\m@ne \fi
   \dg@XTEMP=\dg@DX \dg@changesign\dg@XTEMP\dg@DX
   \dg@YTEMP=\dg@DY \dg@changesign\dg@YTEMP\dg@DY
   \ifnum\dg@YTEMP>\dg@XTEMP \dg@XTEMP=\dg@YTEMP \fi
   \ifnum\dg@XTEMP=\z@ \dg@XTEMP=\@ne \fi
   \multiply\dg@XLBLOFF\dg@LBLOFF \divide\dg@XLBLOFF\dg@XTEMP
   \multiply\dg@YLBLOFF\dg@LBLOFF \divide\dg@YLBLOFF\dg@XTEMP
   \multiply\dg@X\@m \multiply\dg@X\dg@XGRID
   \multiply\dg@Y\@m \multiply\dg@Y\dg@YGRID
   \relax}%

\def\dg@rmcommondiv#1#2#3{%
   \dg@XTEMP=#2\relax
   \divide\dg@XTEMP #1\relax \multiply\dg@XTEMP #1\relax
   \dg@YTEMP=#3\relax
   \divide\dg@YTEMP #1\relax \multiply\dg@YTEMP #1\relax
   \ifnum\dg@XTEMP=#2\relax \ifnum\dg@YTEMP=#3\relax
      \divide#2#1\relax \divide#3#1\relax \fi\fi}%

\def\dg@changesign#1#2{%
   \ifnum #2<\z@ \multiply#1\m@ne
   \else\ifnum #2=\z@ #1=\z@ \fi\fi}%

\def\dg@getnodesize#1#2#3#4#5{%
   #4=\z@\relax #5=\z@\relax
   \expandafter\@for\expandafter\dg@trynode
   \expandafter:\expandafter=#1\do{%
      \dg@XNODE=\m@ne 
      \dg@trynode
      \ifnum #2=\dg@XNODE \ifnum #3=\dg@YNODE
         #4=\dg@XTEMP\relax #5=\dg@YTEMP\relax\fi\fi}}%

\newoptcommand{\dg@makebox}{}[2]{%
   \expandafter\makebox\expandafter(\expandafter
      0\expandafter,\expandafter0\expandafter)\expandafter
      [#1]{#2}}%

\def\dg@novector(#1,#2)#3{}%

\def\dg@letname#1#2{%
   \relax\expandafter
   \let\expandafter #1\csname #2\endcsname\relax}%

\def\dgl@#1{#1{}{}}%
\def\dgl@t#1#2{#1{#2}{}}%
\def\dgl@b#1#2{#1{}{#2}}%
\def\dgl@tb#1#2#3{#1{#2}{#3}}%
\def\dgl@l#1#2{#1{#2}{}}%
\def\dgl@r#1#2{#1{}{#2}}%
\def\dgl@lr#1#2#3{#1{#2}{#3}}%
\makeatother

\theoremstyle{plain}
\newtheorem{thm}{Theorem}[section]
\newtheorem{cor}[thm]{Corollary}
\newtheorem{lem}[thm]{Lemma}

{\theorembodyfont{\rmfamily} \newtheorem{dfn}[thm]{Definition}}
{\theorembodyfont{\rmfamily} \newtheorem{exa}[thm]{Example}}
{\theorembodyfont{\rmfamily} }
{\theorembodyfont{\rmfamily} \newtheorem{que}[thm]{Question}}
{\theorembodyfont{\rmfamily} \newtheorem{rem}[thm]{Remark}}

\DeclareMathOperator{\Grass}{Grass}

\newcommand{\lto}{\longrightarrow}
\newcommand{\V}{\mathbf V}
\newcommand{\G}{\mathbb G}
\newcommand{\R}{\mathcal R}
\newcommand{\Q}{\mathcal Q}
\newcommand{\Gor}{Gor}
\newcommand{\Cat}{Cat}

\newcommand{\Sec}{Sec}
\newcommand{\Ker}{Ker\,}
\newcommand{\rk}{rk\,}

\begin{document}
\title{\bf Chordal varieties of Veronese varieties and
catalecticant matrices }
\date{}
 \author{Vassil Kanev}
\maketitle
\begin{abstract}
It is proved that the  the chordal variety of the Veronese
variety $v_{d}(\mathbb{P}^{n-1})$ is projectively normal,
arithmetically Cohen-Macaulay and its homogeneous ideal  is
generated by the $3\times
3$ minors of two catalecticant matrices. These results are generalized
to the catalecticant varieties $\Gor_{\leq }(T)$ with $t_1=2$.
 \end{abstract}
\section*{Introduction}
  The $r$-secant variety to the Segre variety $\sigma
(\mathbb{P}^{m-1}\times \mathbb{P}^{n-1})$ is  the  projectivization  of
the determinantal variety $M_{r}(m,n)$ of matrices of  rank $\leq r$.
 Its homogeneous
ideal is generated by the $(r+1)\times (r+1)$ minors of the generic $m
\times n$
matrix. Similarly the homogeneous ideal of the $r$-secant variety to
the Veronese variety $v_{2}(\mathbb{P}^{n-1})$  is  generated  by  the
$(r+1)\times (r+1)$ minors of a generic symmetric $n \times n$ matrix.
These two statements are equivalent to the
Second Fundamental Theorem of invariant  theory for the groups $GL_n$
and $O(n)$ (see \cite{Weyl,DP,VP}). In both cases the secant varieties are
projectively normal and arithmetically Cohen-Macaulay (ACM)(see
\cite{BW,Ku}).
\par
The $r$-secant varieties $\Sec _{r}(v_{d}(\mathbb{P}^{n-1}))$ to  the
higher Veronese varieties are related to an  old  problem  of  the
theory of invariants: representing a homogeneous form of degree $d$
in $n$  variables as a sum of $r$  powers of linear forms
\begin{equation}\label{eii}
f = L_{1}^{d}+\cdots +L_{r}^{d}
\end{equation}
Finding   explicitly   polynomial   equations   that    determine
set-theoretically  $\Sec _{r}(v_{d}(\mathbb{P}^{n-1}))$  is  equivalent   to
finding explicit polynomial equations on the  coefficients  of $f$
which are necessary and sufficient for $f$  to  be  representable  in
the form \eqref{eii} or its degeneration. In analogy with the above
examples the following questions arise. What are the generators of
the homogeneous ideal of $\Sec _{r}(v_{d}(\mathbb{P}^{n-1}))$? Do these
varieties satisfy the properties of being projectively normal and ACM?
\par
These questions have affirmative answer in the case of $r$-secants to
a rational normal curve (see Example~\ref{fiv:one} and references
therein). Little is known about these problems when $n\geq 3$. We
answer to them affirmatively in Theorem~\ref{twyeig} for the chordal
variety $\Sec _{2}(v_{d}(\mathbb{P}^{n-1}))$. We prove, assuming the
ground field is algebraically closed of characteristic $0$, that the
homogeneous
ideal of $\Sec _{2}(v_{d}(\mathbb{P}^{n-1}))$ is
generated by the $3\times 3$ minors of the catalecticant
matrices $\Cat _{F}(1,d-1;n)$ and $\Cat _{F}(2,d-2;n)$ (see Section 1 for
definitions); it is projectively normal,
ACM and its affine cone has rational
singularities\footnote{This
result is related to a conjecture posed by A. Geramita in his talk at
the Midwest Algebraic Geometry Conference (Notre Dame, Nov. 7-9,
1997), that the ideal of $\Sec _{2}(v_{d}(\mathbb{P}^{n-1}))$ is
generated by the $3\times 3$ minors of each of the catalecticant matrices
$\Cat _{F}(i,d-i;n),\, 2\leq i\leq d-2$. In the same talk A. Geramita
communicated  a result of T. Deery that these minors generate one and the
same ideal for different $i$ as above when $n = 3,4,5$.}. In
Theorem~\ref{twysix} we generalize these results to some varieties
$\Gor_{\leq }(T)$ which parameterize forms with prescribed dimensions of
the spaces of partial derivatives (see Section 1).
\par
Section 1 contains preliminary material on catalecticant matrices and
their determinantal loci. It is mainly an extract from \cite{IK}.
\par
Section 2 is devoted to the rank $\leq r$ locus of the generic
catalecticant matrix $\Cat _{F}(1,d-1;n)$, which is equivalently the
locus of homogeneous forms of degree $d$ in $n$ variables expressible
in $r$ variables after a linear change of the coordinates. The results
of this section are due to O. Porras \cite{Po} (see also \cite{FW}).
We present here a simplified proof of her theorem and give some
corollaries.
\par
Section 3 contains the results about the chordal variety $\Sec
_{2}(v_{d}(\mathbb{P}^{n-1}))$ that we stated above.
\par
Section 4 contains generalization of the results of Section 3 to the
varieties $\Gor _{\leq }(T)$
 with $t_1 = 2$.
\par
Our approach in proving these results is applying the method
 of G. Kempf \cite{Ke1,Ke2} which was
later developed and used many times for calculating the ideal and the
syzygies of various types of determinantal varieties (see
\cite{La,JPW,We,Po,FW}).
\par
We assume the ground field $k$ is algebraically closed of
characteristic $0$. We mean by variety a separated scheme of finite
type over $k$, i.e. not necessarily irreducible. Unless otherwise
stated by a point of a scheme we mean a closed point of the scheme.

\par
\medskip
\noindent {\bf Acknowledgments.}
The author is grateful to J. Weyman for bringing attention to the
paper of O. Porras \cite{Po} and providing the manuscript \cite{FW},
as well as to A. Iarrobino for useful discussions. The hospitality of
Northeastern University during the author's visit in the Fall of 1997
is gratefully acknowledged.
\section{Catalecticant varieties}\label{sone}
Let $S=k[x_1,\dots ,x_n]$ and let $S_d$ denote the space of
homogeneous polynomials in $S$ of degree $d$. It is a classical
problem of the theory of invariants to find conditions of a form $f\in
S_d$ such that it can be represented as a sum
\begin{equation}\label{eone}
f = L_{1}^{d}+\cdots +L_{r}^{d}
\end{equation}
where $L_1,\dots ,L_d$ are linear forms. Let us denote by $PS(r,d;n)$
the algebraic closure of the set of forms representable as in
\eqref{eone}. When no confusion arises we will use the shorter
notation $P_r$. The Veronese map $v_d : \mathbb{P}(S_1) \to
\mathbb{P}(S_d)$ is given by $v_d(L) = L^d$. We see that $PS(r,d;n)$
is the affine cone of the variety $\Sec _{r}(v_{d}(\mathbb{P}^{n-1}))$
of the $r$-secant $(r-1)$-planes to the Veronese variety
$v_d(\mathbb{P}^{n-1})$. So, finding explicitly polynomial equations
of $\Sec _{r}(v_{d}(\mathbb{P}^{n-1}))$ is equivalent to finding
polynomial equations on the coefficients of $f\in S_d$, which are
satisfied if and only if $f$
has a representation of the form \eqref{eone} or its degeneration (the
later is called generalized additive decomposition \cite{IK}). We
refer the reader to the papers \cite{DK,ER,IK,Ge} for a modern account
and references on the subject.
\par
There is a nice set of polynomials in the ideal of $\Sec
_{r}(v_{d}(\mathbb{P}^{n-1}))$ which we now introduce. Consider the
polynomial ring $R=k[y_1,\dots ,y_n]$ with homogeneous components $R_i$
and consider the differential action of $R$ on $S$ defined as follows.
For $\phi \in R_{d-i}, f\in S_d$ we let
\begin{equation}\label{etwo}
\phi \circ f = \phi (\partial _1,\ldots ,\partial _n)f \in S_i
\end{equation}
Let us choose bases of $R_j$ and $S_j$ as follows. For $R_j$ one takes
\[
Y^V=y_{1}^{v_1}\cdots y_{n}^{v_n}
\]
 with $|V|=v_1+\cdots +v_n=j$ and for $S_j$ one takes
\[
X^{(U)} = \frac{1}{u_{1}!\cdots u_{n}!}x_{1}^{u_1}\cdots x_{n}^{u_n}
\]
with $|U|=u_{1}+\cdots +u_{n}=j$. Fixing a form $f \in S_d$
\[
f = \sum_{|W|=d} a_{W} X^{(W)}
\]
one obtains for every $i, 1\leq i\leq d-1$ a linear map
\begin{equation}\label{ethr}
AP_{f}(i,d-i) : R_{d-i} \lto S_{i},\quad \phi \mapsto \phi \circ f
\end{equation}
which in the bases introduced above has the following matrix
\[
\Cat _{f}(i,d-i;n) = (a_{U,V})_{|U|=i,|V|=d-i},\quad a_{U,V}=a_{U+V}
\]
This is the i-th {\it catalecticant matrix} of $f$. Obviously
\[
^{t}\Cat _{F}(i,d-i;n) = \Cat _{F}(d-i,i;n).
\]
Suppose $f$ has a representation of the form \eqref{eone}. Then for
every $j$ with $1\leq j\leq d-1$ the spaces of $j$-th partial
derivatives $R_j\circ f$ have dimensions $\leq r$. Equivalently the
catalecticant matrices of $f$ have ranks $\leq r$. Let
\[
f = \sum_{|W|=d} Z_{W} X^{(W)}
\]
be a generic form in $S_d$, i.e. the coefficients $Z_{W}$ are
indeterminates over $k$. Then the above can be restated as follows.
The $(r+1)\times (r+1)$ minors of the generic catalecticant matrix
\begin{equation}\label{efou}
\Cat _{F}(i,d-i;n) = (Z_{U+V})_{|U|=i,|V|=d-i}
\end{equation}
vanish on $PS(r,d;n)$ (equivalently on $\Sec
_{r}(v_{d}(\mathbb{P}^{n-1}))$).
\begin{exa}\label{fiv}
Let $d=2$. For a general quadratic form $F =\; ^{t}XZX$ with $X=\;
^{t}(x_1,\ldots ,x_n)$,\break $ Z = (Z_{ij}), \; ^{t}Z = Z$ one has $\Cat
_{F}(1,1;n) = Z$. From linear algebra $P_r=PS(r,2;n) \subset S_2$ is
defined set-theoretically by vanishing of the $(r+1)\times (r+1)$
minors of $(Z_{ij})$. In fact a stronger result is known \cite{Ku}
(see also Corollary~\ref{twy:one}), the ideal $I_r$ of $P_r$ is
generated by the $(r+1)\times (r+1)$ minors of $(Z_{ij})$,
the variety $P_r$ is projectively normal and
arithmetically Cohen-Macaulay.
\end{exa}
\begin{exa}\label{fiv:one}
Let $n=2$. This case was much studied in the 19th century. The
classical references are \cite{GY,El}. A modern exposition can be
found in \cite{KR}. For a general binary form
\[
F = Z_{0}\frac{1}{d!}x_{1}^{d}+\dots
+Z_{j}\frac{1}{(d-j)!j!}x_{1}^{d-j}x_{2}^{j}+\dots
+Z_{d}\frac{1}{d!}x_{2}^{d}
\]
the $i$-th catalecticant matrix is
\[
\Cat _{F}(i,d-i;2) =
\begin{pmatrix}Z_{0}&Z_{1}&\dots &Z_{d-i}\\
Z_{1}&Z_{2}&\dots &Z_{d-i+1}\\
\vdots&\vdots& &\vdots\\
Z_{i}&Z_{i+1}&\dots &Z_{d}
\end{pmatrix}
\]
known also as Hankel matrix. When $P_r=PS(r,d;n)$ is a proper subset
of $S_d$ (this holds iff $2r \leq d$) it is defined set-theoretically
by vanishing of the $(r+1)\times (r+1)$ minors of any of the
catalecticant matrices $\Cat _{F}(i,d-i;2)$ with $r\leq i\leq d-r$
(see e.g. \cite[p. 103]{Ha}). In fact the $(r+1)\times (r+1)$ minors
of any of
the matrices $\Cat _{F}(i,d-i;2)$ with $r\leq i\leq d-r$ generate the
ideal of $P_r$, the variety $P_r$ is projectively normal and
arithmetically Cohen-Macaulay (see \cite{GP,Ei1,Wa}). A detailed
exposition of this material is given also in \cite{IK}.
\end{exa}
We now define two types of catalecticant subvarieties of $S_d$ defined by
vanishing
of certain minors of certain catalecticant matrices.
\begin{dfn}\label{six}
Consider the ideal $J_{r}=I_{r+1}(\Cat _{F}(i,d-i;n))$ generated by the
$(r+1)\times (r+1)$ minors of $\Cat _{F}(i,d-i;n)$ (see \eqref{efou}).
We denote by $\V _{r}(i,d-i;n)$ the closed affine subscheme of $S_d$
whose ideal is $J_r$ and denote by $V_{r}(i,d-i;n)$ the reduced subscheme
\(
\V _{r}(i,d-i;n)_{red}.
\)
\end{dfn}

\begin{dfn}\label{sev}
With every form $f\in S_d$ one associates an ideal $Ann(f)\subset R$
consisting of polynomials $\phi$ such that $\phi \circ f = 0$. These
polynomials are called {\it apolar} to $f$. The algebra $A_f=R/Ann(f)$
is a graded artinian Gorenstein algebra (see e.g. \cite[p. 527]{Ei2} or
\cite{IK}). Its Hilbert sequence is
\[
H(A_{f})=(1,\dots ,h_{i}=\dim _{k}(A_{f})_{i},\dots ,1)
\]
It is symmetric with respect to $\frac{d}{2}$.
\end{dfn}

\begin{dfn}\label{sev:one}
Let $T=(1,t_{1},\dots ,t_{i},\dots ,t_{d-1},1)$ be a sequence
of $d+1$ positive integers with $t_{0}=t_{d}=1$ which is symmetric
with respect to $\frac{d}{2}$. Consider the subset of $S_d$
\[
\Gor (T)\ =\ \{ f\in S_d \; | \; H(A_f)=T\}
\]
The variety $\Gor (T) $ is quiasiaffine being an open subset of
\[
\Gor _{\leq }(T)\ =\ \bigcap_{i=1}^{d-1} V _{t_i}(i,d-i;n)
\]
One puts a finer scheme structure on $\Gor (T) $ by considering the
closed subscheme of $S_d$
\[
\mathbf{Gor}_{\leq}(T)\ =\ \bigcap_{i=1}^{d-1} \V _{t_i}(i,d-i;n)
\]
with ideal $\sum_{i=1}^{d-1}I_{t_i+1}(\Cat _{F}(i,d-i;n))$. Then
$\mathbf{Gor}(T)$ is the open subscheme of $\mathbf{Gor}_{\leq}(T)$
associated with the open set $\Gor (T) $
\end{dfn}
\medskip
Recall we defined a differentiation action $R_i\times S_d \to S_{d-i}$
which is perfect pairing for $i=d$. We will need the following theorem
proved in \cite[Chapter 2]{IK}.
\begin{thm}\label{eig}
Let $f\in S_d$ and let $I=\bigoplus I_i \subset R$ be the graded ideal
of polynomials apolar to $f$.
\par
{\rm (i)}\  Suppose $f\in \V _{r}(i,d-i;n)$. Then we have for the
tangent space at $f$ the equality
\[
T_{f}\V _{r}(i,d-i;n)\ =\ (I_{i}I_{d-i})^{\perp }
\]
\indent {\rm (ii)}\ Suppose $f\in \mathbf{Gor}(T) $. Then
\[
T_{f}\mathbf{Gor}(T)\ =\ (I^{2})_{d}^{\perp }
\]
\end{thm}
The connection between $PS(r,d;n)$ and $\Gor (T) $ is the following.
Let us denote by $r_i=\dim _{k}R_i=\binom{n-1+i}{i}$. Suppose $d=2t$
or $2t+1$ and let $r\leq r_t$. Let $u$ be the maximal number $u\leq t$
such that $r_u\leq r$. As in Definition~\ref{sev:one} consider the
following sequence
\[
T_r\ = \ (1,n,\dots ,r_{u},r,r,\dots ,r,r_u,\dots ,n,1)
\]
Then the above considerations can be reformulated as $PS(r,d;n)\subset
\Gor _{\leq}(T_{r})$ (in fact $\subset \overline{\Gor (T) }$ as shown
e.g. in \cite{IK}). It is proved in \cite{IK} that if furthermore
$r\leq r_{t-1}$, then $PS(r,d;n)$ is an irreducible component of
$\Gor_{\leq}(T_{r})$ and the scheme
$\mathbf{Gor}_{\leq}(T_{r})$ is generically smooth along $PS(r,d;n)$.
This suggests the following two questions.
\begin{que}\label{nin}
For which triples $(r,d,n)$ do the $(r+1)\times (r+1)$ minors of the
catalecticant matrices $\Cat _{F}(i,d-i;n),\; 1\leq i\leq \frac{d}{2}$
determine set-theoretically $PS(r,d;n)$, or equivalently when
$PS(r,d;n)=\Gor _{\leq}(T_{r})$?
\end{que}
\begin{que}\label{ten}
Suppose $r,d,n$ is a triple for which the answer to Question~\ref{nin}
is affirmative. Is it true that the ideal of $PS(r,d;n)$ (equal to the
ideal of $\Sec _{r}(v_{d}(\mathbb{P}^{n-1}))$ is generated by the
$(r+1)\times (r+1)$ minors of the catalecticant matrices $\Cat
_{F}(i,d-i;n)$ for $1\leq i\leq \frac{d}{2}$? In other words when
$PS(r,d;n)=\mathbf{Gor}_{\leq}(T_{r})$?
\end{que}
We have already seen in Examples~\ref{fiv},\ref{fiv:one} that the
answer to these questions is affirmative if $d=2$ or $n=2$. In the
next two sections we will see it is also affirmative if $r\leq 2$.
\par
From representation-theoretic point of view it is convenient to
reformulate some of the notions we have so far encountered in terms of
symmetric tensors. We consider a vector space $V$ of dimension $n$
with a tautological representation of $GL_{n}$. We have
\[
k[x_1,\dots ,x_n]\ =\ S\ \cong \ Sym(V^{*})\ =\ \bigoplus_{d\geq
0}S_{d}V^{*}
\]
Explicitly this isomorphism is: given a homogeneous form $f$ of degree
$d$ there exists a unique symmetric covariant tensor $\tilde{f} \in
S_{d}V^{*}$, such that $f(v)=\tilde{f}(v,\dots ,v) $for every $v\in
V$. Then
\[
R\ \cong \ Sym(V)\ = \ \bigoplus _{i\geq 0}S_{i}V
\]
and the differential action of $R_{d-i}$ on $S_{d}$ defined in
\eqref{etwo} equals $\frac{d!}{i!}$ times the contraction action of
tensors
\[
S_{d-i}V\ \times \ S_{d}V^{*}\ \lto \ S_{i}V^{*}
\]
For later use we also need a coordinate free description of the ideal
generated by the $r\times r$ minors of $\Cat _{F}(i,d-i;n)$. Consider
the linear map
\[
\mu _{r}\ :\ \Lambda ^{r}(S_{d-i}V)\otimes \Lambda ^{r}(S_{i}V)\lto
S_{r}(S_{d}V)
\]
defined as follows:
\[
\mu _{r}((\phi_{1}\wedge \cdots \wedge \phi_{r})\otimes
(\psi_{1}\wedge \cdots \wedge \psi_{r}))\ =\ \det (\phi_{i}\psi_{j}).
\]
Consider the graded ring $Sym(S_{d}V)=\bigoplus_{\ell \geq
0}S_{\ell}(S_{d}V)$. Tensoring by $Sym(S_{d}V)$ one extends $\mu_{r}$
to a homomorphism of degree $0$ of graded $Sym(S_{d}V)$-modules
\begin{equation}\label{eten:b}
\Lambda ^{r}(S_{d-i}V)\otimes \Lambda ^{r}(S_{i}V)\otimes
Sym(S_{d}V)(-r)\lto Sym(S_{d}V)
\end{equation}
\begin{lem}\label{ten:b}
Let us identify as above $S_d$ with $S_dV^{*}$ and the coordinate ring
$k[S_{d}]$ with $Sym(S_{d}V)$. Then
\par
{\rm (i)}\ The image of $\mu_{r}$ is the linear subspace generated by
the $r\times r$ minors of $\Cat _{F}(i,d-i;n)$.
\par
{\rm (ii)}\ The image of the homomorphism \eqref{eten:b} equals the
ideal generated by the $r\times r$ minors of $\Cat _{F}(i,d-i;n)$
\end{lem}
{\bf Proof.}
Let $\{E_{U}:|U|=i\},\; \{E_{W}:|W|=d-i\}$ be the bases of $S_{i}V,\;
S_{d-i}V$ which correspond to the bases $\{Y^{U}\},\; \{Y^{W}\}$ of
$R_{i},\; R_{d-i}$ considered above. Then
\[
\begin{aligned}
&\quad \mu_{r}((E_{W_1}\wedge \cdots \wedge E_{W_r})\otimes
(E_{U_1}\wedge \cdots \wedge E_{U_r}))\\
&= \det (E_{W_i}E_{U_j})\ =\ \det (E_{W_i+U_j})
\end{aligned}
\]
The value of this symmetric tensor on $f\in S_{d}V^{*}$ is equal to
\[
\det (E_{W_i+U_j})\underbrace{(f,\dots ,f)}_{r\text{ times}}\ =\ \det
(\langle E_{W_i+U_j}\: ,\: f\rangle)
\]
The right-hand side is $(\frac{1}{r!})^{r}$ times the $r\times r$
minor of $\Cat _{f}(i,d-i;n)$ corresponding to rows $U_1,\dots ,U_r$
and columns $W_1,\dots ,W_r$. This proves (i). Part (ii) is immediate
from (i).
\hfill $\Box$
\section{Porras' theorem}\label{stwo}
The aim of this section is to give a simplified proof and some
corollaries of a theorem of O. Porras \cite{Po} about the
catalecticant variety $V_{r}(1,d-1;n)$. Porras studies more generally
rank varieties of tensors of arbitrary type. The paper \cite[Sections
4,5]{FW} is focused on the symmetric case and Porras' theorem is
generalized to several symmetric tensors. In \cite[Section 4]{FW} the
reader can find a very clear exposition of the geometric method of
calculating syzygies, which is an important ingredient of the original
proof of Porras' theorem and has many other applications (see
\cite{JPW,PW,We}).
\begin{lem}\label{ele}
Let $1\leq r\leq n-1$.
\par
{\rm (i)}\ Suppose $f\in V_{r}(1,d-1;n)$. Then there is a linear
change of coordinates $x_i=\sum c_{ij}x'_{j}$, such that $f\in
k[x'_{1},\dots ,x'_{r}]_{d}$.
\par
{\rm (ii)}\ $V_{1}(1,d-1;n)\ =\ PS(1,d;n)$
\end{lem}
{\bf Proof.}
(i).\; $\Cat _{f}(d-1,1;n)=\: ^{t}\Cat _{f}(1,d-1;n)$. So, if the rank
of these matrices is $\leq r$, then the kernel of the first one has
dimension $\geq n-r$. Choosing new coordinates $y'_{i}=\sum
c_{ji}y_{j}$ for $R$, so that $y'_{r+1},\dots ,y'_{n}$ belong to this
kernel we have for the dual coordinates of $S,\; x_i = \sum
c_{ij}x'_{j}$ the equation $\partial _{x'_{j}}(f)=0$ for $j\geq r+1$.
This proves (i).
\par
(ii). Immediate from (i).
\hfill $\Box$
\par \medskip
Let $X$ be the affine space $S_{d}V^{*}$ with coordinate ring
$A=k[X]=Sym(S_{d}V)$ with graded components $A_j=S_{j}(S_{d}V)$. The
identification between homogeneous polynomials and symmetric tensors
from the end of Section~\ref{sone} together with Lemma~\ref{ele} yield
that $V_{r}(1,d-1;n)$ equals the variety $X_r$ of symmetric tensors of
rank $r$ \cite{Po,FW}, i.e. tensors which belong to $S_{d}Q^{*}\subset
S_{d}V^{*}$ for some subspace $Q^{*}\subset V^{*}$ of dimension
$r$. In terms of symmetric tensors the catalecticant matrix appears as
follows. A generic tensor $\tilde{F} \in S_{d}V^{*}$ yields by
contraction the map
\[
\Psi \ :\ S_{d-1}V \lto V^{*}\otimes A_{1}
\]
(cf. \cite[p. 703]{Po} and \cite[\S 5]{FW}). Its matrix is up to
constant the catalecticant matrix $\Cat _{F}(1,d-1;n)$.
\begin{thm}[O. Porras]\label{thn}
Let $1\leq r\leq n-1$.
\par
{\rm (i)}\ The variety $X_r=V_{r}(1,d-1;n)$ is irreducible of
dimension $\binom{r+d-1}{d}+r(n-r)$. It is normal, Cohen-Macaulay with
rational singularities.
\par
{\rm (ii)}\ Its ideal $I(X_r)$ equals the ideal $J_r$ generated by the
$(r+1)\times (r+1)$ minors of the catalecticant matrix $\Cat
_{F}(1,d-1;n)$.
\par
{\rm (iii)}\ The singular locus of $X_r$ equals
$X_{r-1}=V_{r-1}(1,d-1;n)$.
\end{thm}
Before giving the proof of this theorem we make some comments on the
original proof (\cite{Po} and \cite{FW}). It consists of three steps.
First one constructs a canonical desingularization $q:Z\to X_r$ and
proves (i) by calculating $R^{i}q_{*}(\mathcal{O}_{Z})$. This step is a
particular
case of a theorem of G. Kempf \cite[p. 239]{Ke2}. Second, using
induction and a
representation-theoretic argument one reduces the proof of the
equality $I(X_{r})=J_r$ to the case $r=n-1$. The third and most
difficult step is to prove $I(X_{n-1})=J_{n-1}$. This is done using a
general theorem \cite[\S 4]{FW} by which one calculates the terms of
the minimal resolution of $I(X_{n-1})$. Our simplification is in the
third step. We obtain easily the result using Theorem~\ref{eig}. For
the sake of completeness we also include the proof of the first two
steps following \cite{FW}.
\par \medskip \noindent
{\bf Proof of Theorem~\ref{thn}.}
\par
 \emph{Step 1}.\: We let $\G =\Grass(n-r,V)=\Grass(r,V^{*}),\;
\mathbb{V}=V\times \G$ and consider the two tautological sequences
\begin{gather*}
0\lto \R \lto \mathbb{V} \lto \Q \lto 0\\
0\lto \Q ^{*}\lto \mathbb{V}^{*} \lto \R ^{*} \lto 0
\end{gather*}
We have $S_{d}\Q^{*}\subset S_{d}\mathbb{V}^{*} =X\times \G$. Let us
denote by $Z$ the total space of the vector bundle $p:S_{d}\Q^{*}\to
\G$. Then by Lemma~\ref{ele} the projection onto the factor $X$ gives
a commutative diagram
\[
\begin{diagram}
\node{Z}\arrow{e,t}{\subset}\arrow{s,l}{q'}\node{X\times
\G}\arrow{s,r}{q}\\
\node{X_r}\arrow{e,t}{\subset}\node{X}
\end{diagram}
\]
with epimorphic $q'$. Thus $X_r$ is irreducible. Furthermore the
restriction of $q'$ on every fiber of the vector bundle
$p:S_{d}\Q^{*}\to \G$ is closed embedding. Hence $q'$ is bijective over
$X_{r}-X_{r-1}$. The later subset is nonempty, since it contains $f =
L_{1}^{d}+\cdots +L_{r}^{d}$ with sufficiently general $L_1,\dots
,L_r$. This proves $Z$ is a resolution of $X_r$, gives the dimension
of $X_r$ and proves $X_{r}- X_{r-1}$ is nonsingular.
\par
The sheaves $R^{i}q'_{*}(\mathcal{O}_{Z})$ are coherent over the
affine scheme $X_r$, so are determined completely by the $A$-module
$H^{i}(Z,\mathcal{O}_{Z})$. The projection $p:Z\to \G$ is affine, thus
\[
H^{i}(Z,\mathcal{O}_{Z})\ =\ H^{i}(\G ,p_{*}\mathcal{O}_{Z}).
\]
Since $Z$ is the total space of the vector bundle $p:S_{d}\Q^{*}\to
\G$ one obtains $p_{*}\mathcal{O}_{Z}=Sym(S_{d}\Q)$. One decomposes
into a direct sum
\[
Sym(S_{d}\Q)\ =\ \bigoplus S_{\lambda }(\Q)
\]
where $S_{\lambda }$ are the Schur functors associated with
certain Young diagrams $\lambda _{0}\geq \cdots \geq \lambda _{r}\geq
0$. Using Bott's theorem (see e.g. \cite[p. 687]{Po} or \cite[p. 232]{We}
for a convenient formulation) one obtains
\[
\begin{aligned}
&H^{i}(\G,p_{*}\mathcal{O}_{Z})\ =\ 0 \text{\quad for\quad}i\geq 1\\
&H^{0}(\G,p_{*}\mathcal{O}_{Z})\ =\ \bigoplus_{\ell(\lambda )\leq
r}S_{\lambda }(V)
\end{aligned}
\]
Here $\ell(\lambda )$ is the number of columns in the Young diagram.
Since $S_{d}V=H^{0}(\G,S_{d}\Q)$ and $k[X]=Sym(S_{d}V)$ we obtain that
$R^{i}q'_{*}\mathcal{O}_{Z}$ is a factor of $k[X]$. This proves $X_r$
is normal with rational singularities which implies by
\cite[p. 50]{KKMS} that it is Cohen-Macaulay.
\par
\emph{Step 2}.\quad Let $J_{r}$ be the ideal generated by the $(r+1)\times
(r+1)$ minors of $\Cat _{F}(1,d-1;n)$. We have $J_r\subset I(X_{r})$
and we want to prove equality. Both ideals are $GL_{n}$-invariant
graded ideals of $k[X]$ and by Step 1 we have
\begin{equation}\label{sin}
I(X_{r})\ =\ \bigoplus_{\ell(\lambda )\geq r+1}S_{\lambda }(V)
\end{equation}
To be more precise this means that in the decomposition of $I(X_{r})$
enter all irreducible components $S_{\lambda }(V)$ of a decomposition
of $k[X]$ with number of columns $\ell(\lambda )\geq r+1$. One proves
$J_r=I(X_{r})$ by descending induction in $r$, the case $r=n$ being
trivial. Suppose one has that $J_{r+1}=I(X_{r+1})$. Using
\eqref{sin} one has to prove that every isotypical component of
$k[X]$ associated to a $S_{\lambda }(V)$ with $\ell(\lambda )=r$ is
contained in $J_r$. We claim it suffices to prove the equality
$J_r=I(X_{r})$ in the case of vector spaces $V$ of dimension $r+1$.
Indeed, suppose $r+1<n$. Choose a basis $e_1,\dots,e_n$ of $V$, let
$V'=\langle e_1,\dots,e_{n-1}\rangle$, let $GL_{n-1}\subset GL_{n}$ be
the corresponding embedding and let $V\to V'$ be the $GL_{n-1}$-invariant
projection. It induces a projection
\[
k[X]\ =\ Sym(S_{d}V)\lto Sym(S_{d}V')\ =\ k[X']
\]
with the property that every isotypical component
$\oplus_{1}^{m}S_{\lambda }(V)$ of $k[X]$ with $\ell(\lambda )\leq
n-1$ is transformed into isotypical component
$\oplus_{1}^{m}S_{\lambda }(V')$ of $k[X']$. For homogeneous forms in
$S_d\cong S_{d}V^{*}$ the effect of this projection is letting
$x_n=0$, so the images of $J_r$ and $I(X_{r})$ are the corresponding
ideals for polynomials in $n-1$ variables. Thus, if $r+1\leq n-1$
proving the equality
\[
I_{r+1}(\Cat _{F'}(1,d-1;n-1))\ =\ I(V_{r}(1,d-1;n-1))
\]
would imply that every isotypical component $\oplus_{1}^{m}S_{\lambda
}(V)$ with $\ell(\lambda )=r+1$ is contained in $J_r$ which would
imply from the induction hypothesis $J_{r+1}=I(X_{r+1})$ that
$J_{r}=I(X_{r})$. Repeating this argument we see it suffices to prove
that $J_{r}=I(X_{r})$ for $n=r+1$.
\par
\emph{Step 3}.\quad Let $r=n-1$. Then the ideal $J_{n-1}$ is generated by the
maximal minors of the catalecticant matrix $\Cat _{F}(1,d-1;n)$ of
type $n\times N$ with $N=\binom{n+d-2}{d-1}$. The codimension of
$X_{n-1}$ equals
\[
\binom{n+d-1}{d}-\binom{n+d-2}{d}-(n-1)\ =\ N-n+1.
\]
This number is the codimension of the generic determinantal locus
$M_{n-1}(n,N)$, hence the catalecticant scheme $\V _{r}(1,d-1;n)$ with
ideal $J_r$ is Cohen-Macaulay (see e.g. \cite[p. 84]{ACGH}). Since $\V
_{r}(1,d-1;n)$ is an irreducible scheme, in order to prove that $J_r$
is a prime ideal (and thus equal to $I(X_{r})$) it suffices to verify
that $\V _{r}(1,d-1;n)$ is generically smooth (see e.g.
\cite[p. 457]{Ei2}). We prove smoothness at every $f\in \V
_{n-1}(1,d-1;n)- \V _{n-2}(1,d-1;n)$. Let $I\subset R$ be the
ideal of polynomials  apolar  to $f$. Then $\dim _{k}I_1=1,\: \dim
_{k}I_{d-1}=\dim _{k}R_{d-1}-n+1$ since $\rk \Cat _{f}(1,d-1;n)=n-1$.
Thus
\[
\dim _{k}I_{1}I_{d-1}\ =\ \binom{n+d-2}{d-1}-n+1
\]
which by Theorem~\ref{eig} yields that $f$ is a smooth point of the
scheme $\V _{n-1}(1,d-1;n)$. This proves (ii).
\par
It remains to prove that $Sing (X_r) = X_{r-1}$. We have already proved
above that $X_r- X_{r-1}$ is nonsingular. Associating to
$f\in S_d$ the catalecticant matrix $\Cat _{f}(1,d-1;n) $ yields a
linear map into the space of $n\times N$ matrices
$Cat : S_d \lto M(n,N)$. The $(r+1)\times (r+1)$ minors are
homogeneous polynomials on $M(n,N)$ which vanish of order $\geq 2$ at
every point of the rank $\leq r-1$ locus $M_{r-1}(n,N)$ (see
\cite[pp. 184-185]{Ha}). The same holds for the pull-back of these
minors by the map $Cat$, so from the equality $J_r=I(X_{r})$ we
conclude that for every $f\in X_{r-1}$ one has $T_fX_r=S_d$. Therefore
$Sing X_r = X_{r-1}$. Theorem~\ref{thn} is proved.
\hfill $\Box$
\begin{rem}\label{nen:a}
Porras' theorem gives affirmative answer to the first two questions of
Problem 11.6 in \cite{Ge}.
\end{rem}
A particular case of Porras' theorem is the following well-known
theorem mentioned in the introduction.
\begin{cor}\label{twy:one}
Let $X$ be the affine space of symmetric $n\times n$ matrices. Then
the sublocus of matrices of rank $\leq r$ is irreducible, normal,
Cohen-Macaulay with rational singularities and its ideal is generated
by the $(r+1)\times (r+1)$ minors of a generic symmetric $n\times n$
matrix. The singular locus of $X_r$ equals $X_{r-1}$.
\end{cor}
\begin{cor}\label{twy:two}
The Veronese variety $v_d(\mathbb{P}^{n-1})$ is projectively normal,
arithmetically Cohen-Macaulay, its affine cone has rational
singularities and its ideal is generated by the $2\times 2$ minors of
the catalecticant matrix $\Cat _{F}(1,d-1;n)$.
\end{cor}
\begin{rem}\label{twy:thr}
The generation of the ideal of $v_d(\mathbb{P}^{n-1})$ by the $2\times
2$ minors was posed as a question in \cite[pp. 101-102]{Ge}. It is
reported in \cite[p. 213]{GPS} that M. Pucci has shown independently
this fact.
\end{rem}
\begin{cor}\label{twy:fou}
The chordal variety $\Sec _{2}(v_{3}(\mathbb{P}^{n-1}))$ is
projectively normal, arithmetically Cohen-Macaulay, its affine cone
has rational singularities and its ideal is generated by the $3\times
3$ minors of the catalecticant matrix $\Cat_{F}(1,2;n)$. Its singular
locus equals the Veronese variety $v_3(\mathbb{P}^{n-1})$.
\end{cor}
{\bf Proof.} It is
proved in \cite{IK} that $\dim _{k}PS(2,d;n)=2n$ which equals the
dimension of $V_{2}(1,2;n)$ according to Theorem~\ref{thn}(i).
\hfill $\Box$
\begin{rem}
This corollary gives affirmative answer to Problem 10.7 in
\cite[p. 102]{Ge}.
\end{rem}
\section{The chordal variety to a Veronese variety}\label{sthr}
The chordal variety to $v_3(\mathbb{P}^{n-1})$ was considered in the
previous section. We assume in the present one $d\geq 4$. We saw in
Section~\ref{sone} that the affine cone to the chordal variety $\Sec
_{2}(v_{d}(\mathbb{P}^{n-1}))$ which is $PS(2,d;n)$ is contained in
$\Gor _{\leq}(T_{2})$, the latter variety being defined by the
vanishing of the $3\times 3$ minors of the catalecticant matrices
$\Cat _{F}(i,d-i;n)$ for $1\leq i\leq \frac{d}{2}$. In fact a smaller
set of equations suffices (see \cite[p. 107]{Ge} for the case $d=3$).
\begin{lem}\label{twytwo}
A form $f\in S_d$ belongs to $PS(2,d;n)=P_2$ if and only if
\[
\rk \Cat _{f}(1,d-1;n)\ \leq \ 2 \text{\quad and\quad } \rk \Cat
_{f}(2,d-2;n)\ \leq \ 2.
\]
The possibilities for $f\in P_2$ are:
$0,L^{d},L_{1}^{d}+L_{2}^{d}\text{\: or\: } L_{1}L_{2}^{d-1}$ for some
linear forms $L,L_{1},L_{2}$.
\end{lem}
{\bf Proof.}
By Lemma~\ref{ele} the first condition gives that after a
change of coordinates $f\in k[x'_{1},x'_{2}]_{d}$. The second one holds
 if and only if there is a form $\phi \in k[y'_1,y'_2]_2$ such that
$\phi \circ f = 0$. For every $s\leq \frac{d}{2}$ the variety
$V_{s}(s,d-s;2)$ is irreducible. Indeed, one constructs its smooth
resolution as follows. One considers
\begin{equation}\label{etwythr:one}
\tilde{Y}\ \subset \ \mathbb{P} (k[y'_1,y'_2]_s \times k[x'_1,x'_2]_d)
\end{equation}
consisting of pairs $([\phi],f)$ with $\phi \circ f = 0$. Then the
first and second projection yield a diagram
\begin{equation}\label{etwythr:two}
\begin{diagram}
\node{\tilde{Y}} \arrow{e,t}{\pi _{2}} \arrow{s,l}{\pi
_{1}}\node{V_{s}(s,d-s;2)}\\
\node{\mathbb{P} (k[y'_1,y'_2]_s)}
\end{diagram}
\end{equation}
such that $\pi _{1}$ is a vector bundle and $\pi _{2}$ is a birational
morphism. The varieties $V_{s}(s,d-s;2)$ have open dense subsets
consisting of forms $f = L_{1}^{d}+\cdots +L_{s}^{d}$ with
nonproportional $L_i\in k[x'_1,x'_2]_1$, these forms being apolar to
$\phi \in k[y'_1,y'_2]_s$ with $s$ distinct roots (see e.g.
\cite{GY,El,IK}). In our particular case $s=2$ we obtain $f\in
PS(2,d;n)$ and if $f\neq 0$ one has the following cases
\par
(a)\: $f=L^{d}$ \; if\; $\rk \Cat _{f}(1,d-1;n)=1$
\par \noindent
and if \; $\rk \Cat _{f}(1,d-1;n)=2$
\par
(b)\: $f=L_{1}^{d}+L_{2}^{d}$ if $\phi $ has $2$ distinct roots;
\par
(c)\: $f=L_{1}L_{2}^{d-1}$ if $\phi $ has a double root.
\par \noindent This proves the lemma.
\hfill $\Box$ \par \medskip \noindent
Let us make the relative analog of the above. We
consider
the diagram as in the proof of Theorem~\ref{thn}
\begin{equation}\label{etwyfou}
\begin{diagram}
\node{Z}\arrow{e,t}{p}\arrow{s,l}{q'}\node{\G}\\
\node{X_2}\arrow{e,t}{\subset }\node{X}
\end{diagram}
\end{equation}
and the map of vector bundles over $Z$
\begin{equation}\label{etwyfiv}
\varphi : S_{d-2}(p^{*}\Q )\lto S_2(p^{*}\Q ^{*})
\end{equation}
which is defined on the fiber over $f\in S_d(\Q ^{*})=Z$ as
contraction with $f$. If we let $A=S_{d-2}(p^{*}\Q ), E=S_2(p^{*}\Q
^{*})$ with ranks $a=d-1, e=3$ respectively we are in the situation
considered in \cite[p. 216]{La}. Let us denote by $Y\subset Z$ the
corank 1 determinantal subscheme of the map $\varphi $, defined as the
closed subscheme of $Z$ with ideal sheaf $J_Y$ generated by the
$3\times 3$ minors of $\varphi $, i.e. $J_Y$ is the image of
\[
\varphi ^{\sharp} : \Lambda ^{3}A \otimes \Lambda ^{3}E^{*} \lto
\mathcal{O}_{Z}
\]
\begin{thm}\label{twysix}
The scheme $Y$ is integral, normal, Cohen-Macaulay with rational
singularities. The resolution of $\mathcal{O}_{Y}$ is given by the
Eagon-Northcott complex.
\begin{multline*}
0\to \Lambda ^{a}A\otimes S_{a-e}E^{*}\otimes \Lambda ^{e}E^{*}\to
\cdots \to \Lambda ^{j}A\otimes S_{j-e}E^{*}\otimes \Lambda ^{e}E^{*}\to
\\
\cdots \to \Lambda ^{e}A\otimes \Lambda ^{e}E^{*}\to \mathcal{O}_{Z}\to
\mathcal{O}_{Y}\to 0
\end{multline*}
where $d-1=a\geq j\geq e=3,\, A=p^{*}S_{d-2}\Q,\, E^{*}=p^{*}S_{2}\Q$.
\end{thm}
{\bf Proof.} The Grassman bundle $G_1(E)$ of quotients of rank $1$
considered in \cite[p. 216]{La} is by duality equal to
\[
\Grass (1,E^{*}) = \mathbb{P} (p^{*}S_2\Q ) = \mathbb{P} (S_2\Q
)\times _{\G} Z
\]
The relative (over $\G $ ) analog of \eqref{etwythr:one},
\eqref{etwythr:two} is a variety
\begin{equation}\label{etwysix:one}
\tilde{Y}\ \subset \ \mathbb{P}(S_2\Q )\ \times _{\G}\ Z
\end{equation}
with two projections $\pi _{1},\pi = \pi _{2}$
\begin{equation}\label{etwysix:two}
\begin{diagram}
\node{\tilde{Y}}\arrow{e,t}{\pi }\arrow{s,l}{\pi _{1}}\node{Y}\\
\node{\mathbb{P}(S_2\Q )}\arrow{e}\node{\G}
\end{diagram}
\end{equation}
The projection $\pi _{1}$ makes $\tilde{Y}$ a vector bundle over
$\mathbb{P}(S_2\Q)$, so it is smooth. The second one $\pi _{2}=\pi$
maps $\tilde{Y}$ birationally onto $Y$, so all condition required in
Proposition 2.4 and Theorem 2.6 of \cite{La} are satisfied. The
arguments in \cite[p. 217]{La} and \cite[pp. 181-182]{Ke1} yield
$R^{p}\pi _{*}\mathcal{O}_{\tilde{Y}}=0$ for $p\geq 1,\: \pi
_{*}\mathcal{O}_{\tilde{Y}}=\mathcal{O}_{Y}$ and easy calculation of
the terms $\mathcal{E}_{1}^{-j,i}$ of the spectral sequence considered in
\cite[p. 217]{La} yield the Eagon-Northcott complex. The remaining
statements of the theorem follow from \cite[p. 50]{KKMS}.
\hfill $\Box$ \par \medskip \noindent
Let $P_2=PS(2,d;n)$. We obtain from \eqref{etwyfou} the following
diagram
\begin{equation}\label{etwysev}
\begin{diagram}
\node{Y}\arrow{e,t}{\subset
}\arrow{s,l}{q}\node{Z}\arrow{e,t}{p}\arrow{s,r}{q'}\node{\G}\\
\node{P_2}\arrow{e,t}{\subset }\node{X_2}\arrow{e,t}{\subset }\node{X}
\end{diagram}
\end{equation}
\begin{thm}\label{twyeig}
Let $P_2=PS(2,d;n)\subset S_d$.
\par
{\rm (i)}\ The variety $P_2$ is normal, Cohen-Macaulay with rational
singularities.
\par
{\rm (ii)}\ The ideal $I(P_2)\subset k[S_d]$ is generated by the
$3\times 3$ minors of the catalecticant matrices $\Cat _{F}(1,d-1;n)$
and $\Cat _{F}(2,d-2;n)$.
\par
{\rm (iii)}\ The singular locus of $P_2$ equals $P_1=PS(1,d;n)$.
\end{thm}
{\bf Proof.} The case $d=3$ is proved in Corollary~\ref{twy:fou}. So
we may assume $d\geq 4$. The map $q':Z\to X_2$ is biregular over
$X_2-X_1$ and $X_1=V_{1}(1,d-1;n)$ is properly contained in $P_2$
(Lemma~\ref{ele}(ii)). Thus $q:Y\to P_2$ is birational. By
\cite[Theorem 5]{Ke2} in order to verify (i) it suffices to prove
that $R^{i}q_{*}\mathcal{O}_{Y}=0$ for $i\geq 1$ and
$q_{*}\mathcal{O}_{Y}=P_2$. Since $P_2$ is a closed affine subvariety of
$X$ it is equivalent to prove that $H^{i}(Y,\mathcal{O}_{Y})=0$ for
$i\geq 1$ and $H^{0}(Y,\mathcal{O}_{Y})$ is a factor of
$k[X]=Sym(S_{d}V)$. The map $p:Z\to \G$ is affine, so pushing
forward the Eagon-Northcott resolution of $\mathcal{O}_{Y}$ one
obtains a resolution for $p_{*}\mathcal{O}_{Y}$:
\begin{multline}\label{etwyeig}
\cdots \to \Lambda ^{j}(S_{d-2}\Q)\otimes S_{j-3}(S_{2}\Q)\otimes
\Lambda ^{d-1}(S_{2}\Q)\otimes Sym(S_{d}\Q)(-j)\to\\
\cdots \to \Lambda ^{3}(S_{d-2}\Q)\otimes \Lambda ^{3}(S_{2}\Q)\otimes
Sym(S_{d}\Q)(-3)\to Sym(S_{d}\Q)\to p_{*}\mathcal{O}_{Y}\to 0
\end{multline}
where $d-1\geq j\geq 3$. With the grading $Sym(S_{d}\Q)(-j)$ the
differentials of the complex are of degree $0$. We need a lemma.
\begin{lem}\label{twynin}
Let $\mathcal{F}_{j}=\bigoplus_{n}(\mathcal{F}_{j})_{n}$ be the $j$-th
graded sheaf of the resolution of $p_{*}\mathcal{O}_{Y}$. Then
$H^{i}(\G ,\mathcal{F}_{j})=0$ for $i\geq 1$.
\end{lem}
{\bf Proof.} This is immediate from Bott's theorem since each sheaf
$(\mathcal{F}_{j})_{n}$ decomposes as direct sum of $S_{\lambda }(\Q)$
for some Young diagrams $\lambda _{0}\geq
\lambda _{1}
\dots \geq\lambda _{r}\geq 0$.
\hfill $\Box$ \par \medskip \noindent
{\it Proof continued.} Splitting the complex $\mathcal{F}_{\bullet }$
into a sequence of short exact sequences we deduce from the lemma
using induction that the higher cohomology groups of each syzygy sheaf
is $0$, hence
\[
H^{i}(Y,\mathcal{O}_{Y})\ =\ H^{i}(\G ,p_{*}\mathcal{O}_{Y})\ =\ 0
\text{\: for\: } i\geq 1
\]
and there is an exact sequence
\begin{equation}\label{ethy}
0\to \Gamma (\mathcal{F}_{d-3})\to \cdots \to \Gamma
(\mathcal{F}_{1})\to \Gamma (p_{*}\mathcal{O}_{Z})\to \Gamma
(p_{*}\mathcal{O}_{Y})\to 0
\end{equation}
One obtains that $H^{0}(Y,\mathcal{O}_{Y})$ is a factor of
$H^{0}(Z,\mathcal{O}_{Z})$ which equals $H^{0}(X_2,\mathcal{O}_{X_2})$
by Theorem~\ref{thn}. This proves (i).
\par
Borel-Weyl's theorem yields the following diagram of maps
\[
\begin{diagram}
\node{\Lambda ^{3}(S_{d-2}V)\otimes \Lambda ^{3}(S_{2}V)\otimes
Sym(S_{d}V)(-3)}\arrow{e,t}{\mu }\arrow{s}\node{Sym(S_{d}V)}\arrow{s}\\
\node{\Gamma (\Lambda ^{3}(S_{d-2}\Q)\otimes \Lambda
^{3}(S_{2}\Q)\otimes Sym(S_{d}\Q)(-3))}\arrow{e,t}{d_{1}}\node{\Gamma
(Sym(S_{d}\Q))}
\end{diagram}
\]
The map $\mu $ is the one defined in \eqref{eten:b}, the map $d_1$ is
a differential of the complex \eqref{ethy} and its image is $\Gamma
(p_{*}J_{Y})$. The vertical maps are epimorphisms induced from the
isomorphism $V\overset{\sim }\lto \Gamma (\G ,\Q)$. The diagram is
commutative as one easily sees by restricting the sections of the
bottom line to an arbitrary element of the Grassmanian $\G$. This
implies by Lemma~\ref{ten:b} that $\Gamma (p_{*}\mathcal{O}_{Y})$ is
generated by the images of the $3\times 3$ minors of $\Cat
_{F}(2,d-2;n)$. Now, $\Gamma (\G ,p_{*}J_{Y})=\Gamma (Z,J_{Y})$ and
the isomorphism $q^{\prime *}:\Gamma
(X_2,\mathcal{O}_{X_2})\overset{\sim }\lto \Gamma (Z,\mathcal{O}_{Z})$
transforms $\Gamma (X_2,J_{P_2})$ onto $\Gamma (Z,J_{Y})$. We conclude
that the images of the $3\times 3$ minors of $\Cat _{F}(2,d-2;n)$ in
$k[X_2]$ generate the ideal $\Gamma (X_2,J_{P_2})$, hence these minors
together with the $3\times 3$ minors of $\Cat _{F}(1,d-1;n)$ generate
the ideal $I(P_2)\subset k[X] = k[S_d]$ (Theorem~\ref{thn}). This
proves (ii).

For the proof of (iii) we refer to Theorem~\ref{thyfiv}(iii).
\hfill $\Box$

\section{The varieties $\Gor _{\leq}(T)$}\label{sfou}
The varieties $PS(2,d;n)$ studied in Section~\ref{sthr} are particular
cases of the varieties $\Gor _{\leq}(T)$ with $t_1=2$ (see
Definition~\ref{sev:one}).
\begin{lem}\label{thythr}
Suppose $f\in V_{2}(1,d-1;n)-V_{1}(1,d-1;n)$. Then there is a number
$s\geq 2$ with $2s\leq d+2$ such that the Hilbert sequence of
$A_f=R/Ann(f)$ is
\begin{equation}\label{ethythr}
H(A_f)\ =\ T_{2,s}\ =\ (1,2,\ldots ,\overset{s-1}s,s,\ldots
,\overset{d-s+1}s,s-1,\ldots ,2,1)
\end{equation}
\end{lem}
{\bf Proof.} By Lemma~\ref{ele} after a change of variables $f$
becomes a binary form and the above Hilbert sequence is one of the
possible Hilbert sequences for binary forms (see e.g. \cite{IK}).
 \hfill $\Box$ \par \medskip \noindent
We see that the only possible Hilbert sequences $T$ with $t_1=2$ are
$T_{2,s}$ defined above. The forms $f$ which have $H(A_f)=T_{2,s}$ can
be explicitly described as follows. Let $z_1,z_2$ be a basis of the
space of linear forms that vanish on $\Ker \Cat _{f}(d-1,1;n)$. Then
there exist linear forms $L_i=\alpha _iz_1+\beta _iz_2,\, i=1,\dots
,m$ not proportional to each other, and polynomials $G_i(z_1,z_2)$
not divisible by $L_i$, of degrees $d_i-1$  such that
\[
f = G_1L_{1}^{d-d_1+1}+\cdots +G_mL_{m}^{d-d_m+1}
\]
and $s=\sum d_i$. Furthermore this representation is unique up to
multiplication of the $L_i$'s by nonzero constants if $2s\leq d+1$
(see \cite{GY} or \cite{IK}).
\begin{lem}\label{thyfou}
Let $s\geq 2,\, 2s\leq d+2$. A form $f\in S_d$ belongs to $\Gor
_{\leq}(T_{2,s})$ if and only if
\[
\rk \Cat _{f}(1,d-1;n)\leq 2,\quad \rk \Cat _{f}(s,d-s;n)\leq s
\]
The variety $\Gor _{\leq}(T_{2,s})$ is irreducible.
\end{lem}
{\bf Proof.} The proof is left to the reader. The same argument as
that of Lemma~\ref{ele} works using \eqref{etwythr:one} and
\eqref{etwythr:two}.
\hfill $\Box$ \par \medskip \noindent
Notice that for the maximal possible values of $s$, namely
$s=\frac{d+1}{2}$ if $d$ is odd and $s=\frac{d+2}{2}$ if $d$ is even,
the rank condition $\rk \Cat _{f}(s,d-s;n)\leq s$ is fulfilled
automatically, so in these cases $\Gor
_{\leq}(T_{2,s})=V_{2}(1,d-1;n)$, the variety studied in
Section~\ref{stwo}. The following theorem generalizes
Theorem~\ref{twyeig}.
\begin{thm}\label{thyfiv}
Let $s\geq 2,\, 2s\leq d$. Consider $\Gor _{\leq}(T_{2,s})$ (see
\eqref{ethythr} and Definition~\ref{sev:one}). Then

{\rm (i)}\ The variety $\Gor _{\leq}(T_{2,s})$ is normal,
Cohen-Macaulay with rational singularities.

{\rm (ii)}\ The ideal of $\Gor _{\leq}(T_{2,s})$ is generated by the
$3\times 3$ minors of $\Cat _{F}(1,d-1;n)$ and the $(s+1)\times (s+1)$
minors of $\Cat _{F}(s,d-s;n)$.

{\rm (iii)}\ If $s\geq 3$ the singular locus of $\Gor
_{\leq}(T_{2,s})$ equals $\Gor _{\leq}(T_{2,s-1})$. The singular locus
of $\Gor _{\leq}(T_{2,2})=PS(2,d;n)$ equals $PS(1,d;n)$.
\end{thm}
{\bf Proof.} The proof of (i) and (ii) is essentially the same as that
of Theorem~\ref{twyeig} and we leave the details to the reader. We
only indicate what changes one needs to do. Instead of \eqref{etwyfiv}
one considers the contraction map
\[
\varphi : S_{d-s}(p^{*}\Q)\lto S_{s}(p^{*}\Q ^{*})
\]
and the rank-$s$ determinantal subscheme $Y\subset Z$ whose ideal
sheaf is generated by the $(s+1)\times (s+1)$ minors of $\varphi $,
i.e. by the image of
\[
\varphi ^{\sharp} : \Lambda ^{s}(S_{d-i}A)\otimes \Lambda
^{s}(S_{d-i}E^{*})\lto \mathcal{O}_{Z}
\]
where $A=S_{d-s}(p^{*}\Q)$ has rank $a=d-s+1$ and $E=S_s(p^{*}\Q ^{*}$
has rank $e=s+1$. The same statements for $Y$ as in
Theorem~\ref{twysix} hold with the Eagon-Northcott complex defined by
$A,E,a$ and $e$ as above.

Next, one replaces $P_2$ by $\Gor _{\leq}(T_{2,s})$ in the diagram
\eqref{etwysev} and pushes forward by $p:Z\to \G$ the Eagon-Northcott
complex. Then exactly by the same arguments as in the proof of
Theorem~\ref{twyeig} one verifies(i) and (ii).

For the proof of (iii) one observes that the open subset $U \subset
\Gor _{\leq}(T_{2,s})$ consisting of $f$ with $\rk \Cat
_{f}(2,d-2;n)=s$ is nonsingular. This follows from the fact that the
birational projection $Y\to \Gor _{\leq}(T_{2,s})$ is biregular over
$U$ and the preimage of $U$ is nonsingular in $Y$ as evident from the
rank-$s$ analog of the diagram \eqref{etwysix:two}. For the proof that
$\Gor _{\leq}(T_{2,s-1})$ (resp. $PS(1,d;n)$ for $s=2$) belongs to the
singular locus of $\Gor _{\leq}(T_{2,s})$ one uses the same argument
as that of Theorem~\ref{thn}(iii) applied to the map
\[
Cat : X_2 \lto M\left( \binom{n+1}{2},\binom{n+d-3}{d-2}\right)
\]
given by the catalecticant matrix $\Cat _{F}(2,d-2;n)$.
\hfill $\Box$

\bibliographystyle{amsalpha}

\bigskip

\noindent
{\sc
Institute of Mathematics, Bulgarian Academy of Sciences,\\
Acad. G. Bonchev Str. bl. 8, Sofia, Bulgaria 1113}

{\it E-mail address:} {\bf kanev@@math.acad.bg}

\end{document}